\documentclass[a4paper, 12pt]{amsart}

\usepackage[latin1]{inputenc}
\usepackage{amsmath}
\usepackage{amsfonts}
\usepackage{amssymb}
\usepackage{amscd} 
\usepackage{mathrsfs}
\usepackage[english]{babel}
\usepackage[size=2em,nohug,midshaft,scriptlabels,PostScript=dvips]{diagrams}

\theoremstyle{plain}
\newtheorem{theorem}{Theorem}[section]
\newtheorem{lemma}[theorem]{Lemma}
\newtheorem{proposition}[theorem]{Proposition}
\newtheorem{corollary}[theorem]{Corollary}

\theoremstyle{definition}
\newtheorem{definition}[theorem]{Definition}
\newtheorem{remark}[theorem]{Remark}

\newcommand{\ZZ}{\mathbb{Z}} 
\newcommand{\QQ}{\mathbb{Q}} 
\newcommand{\iso}{\cong}      
\newcommand{\PP}{\mathbb{P}}  
\renewcommand{\AA}{\mathbb{A}}  
\newcommand{\sheaf}[1]{\mathscr{#1}} 
\newcommand{\OO}{\sheaf{O}}   
\newcommand{\floor}[1]{\lfloor{#1}\rfloor} 
\newcommand{\ceil}[1]{\lceil{#1}\rceil}    
\newcommand{\dual}[1]{{#1}^\vee} 
\newcommand{\tsum}{{\textstyle\sum}}  
\newcommand{\res}[2]{\left.#1\right|_{#2}} 
\newcommand{\isoto}{\overset{\sim}{\to}} 
\newcommand{\abuts}{\Rightarrow} 
\newcommand{\Gm}{\mathbb{G}_m} 

\newcommand{\spot}{{\scriptscriptstyle\bullet}}
\newcommand{\alt}{{\textstyle\wedge}}
\newcommand{\Alt}{{\text{\small $\textstyle\bigwedge$}}}

\newcommand{\union}{\cup}
\newcommand{\Union}{\textstyle{\bigcup}}
\newcommand{\intersect}{\cap}

\newcommand{\tensor}{\otimes}

\newcommand{\plus}{\oplus} 
\newcommand{\Plus}{\textstyle{\bigoplus}}

\DeclareMathOperator{\Pic}{Pic} 
\DeclareMathOperator{\NS}{NS} 
\DeclareMathOperator{\Hom}{Hom}

\DeclareMathOperator{\Ext}{Ext}
\DeclareMathOperator{\Spec}{Spec}

\DeclareMathOperator{\Ker}{Ker}
\DeclareMathOperator{\Coker}{Coker}
\DeclareMathOperator{\Image}{Im}
\DeclareMathOperator{\id}{id}
\DeclareMathOperator{\ob}{ob} 
\newcommand{\cat}[1]{\mathcal{#1}} 
\newcommand{\into}{\hookrightarrow} 
\newcommand{\onto}{\twoheadrightarrow} 

\newcommand{\koszA}{\left(\begin{smallmatrix}
\vartheta^+_i\\
\vartheta^-_i
\end{smallmatrix}\right)}
\newcommand{\koszB}{\left(\begin{smallmatrix}
\vartheta^-_i & -\vartheta^+_i
\end{smallmatrix}\right)}

\author{Martin G. Gulbrandsen}
\address{Stord/Haugesund University College, Norway}
\email{martin.gulbrandsen@hsh.no}

\title{Vector bundles and monads on abelian threefolds}

\begin{document}

\thanks{The author thanks the Max-Planck-Institut f\"ur Mathematik in Bonn
for its hospitality and financial support. The bulk of this paper was written
during a stay as postdoc at the MPIM in 2008--2009.}

\begin{abstract}
Using the Serre construction, we give examples of
stable rank $2$ vector bundles on
principally polarized abelian threefolds $(X,\Theta)$ with Picard number $1$. The
Chern classes $(c_1, c_2)$ of these examples realize roughly one half of the
classes that are a priori allowed by the Bogomolov inequality and Riemann-Roch
(the latter gives a certain divisibility condition).

In the case of even $c_1$, we study deformations of these vector bundles $\sheaf{E}$,
using
a second description
in terms of monads, similar to the ones studied by Barth--Hulek on projective
space. By an explicit analysis of the hyperext spectral sequence
associated to the monad, we show that the space of first order infinitesimal
deformations of $\sheaf{E}$ equals the space of first order infinitesimal
deformations of the monad. This leads to the formula
\begin{equation*}
\dim \Ext^1(\sheaf{E},\sheaf{E}) = \tfrac{1}{3} \Delta(\sheaf{E})\cdot \Theta + 5
\end{equation*}
(we emphasize that its validity is only proved for special bundles $\sheaf{E}$
coming from the Serre construction),
where $\Delta$ denotes the discriminant $4c_2 - c_1^2$.

Finally we show that, in the first nontrivial example of the above construction
(where $c_1=0$ and $c_2 = \Theta^2$), the infinitesimal identification
between deformations of $\sheaf{E}$ and of the monad can be extended to a Zariski
local identification: this leads to an explicit description of a Zariski open 
neighbourhood of $\sheaf{E}$ in its moduli space $M(0,\Theta^2)$. This neighbourhood
is a ruled, nonsingular variety of dimension $13$, birational to a $\PP^1$-bundle
over a finite quotient of $X^2\times_X X^2\times_X X^2$, where $X^2$ is considered
as a variety over $X$ via the group law.
\end{abstract}

\maketitle

\setcounter{tocdepth}{1}
\tableofcontents

\section{Introduction}

The geometry of moduli spaces for stable vector bundles on Calabi-Yau (by which
we just mean having trivial canonical bundle) threefolds is largely unknown,
but of high interest, for instance due to their relevance for string theory,
and the computation of Donaldson-Thomas invariants. In lower dimension, vector
bundles on Calabi-Yau curves (i.e.~elliptic) were classified by Atiyah, and
are parametrized by the same curve. Moduli spaces for stable bundles, or coherent
sheaves, on Calabi-Yau surfaces (i.e.~K3 or abelian) are holomorphic symplectic
varieties (Mukai \cite{mukai84}, generalizing Beauville \cite{beauville83},
generalizing Fujiki \cite{fujiki83}). This is a very rare geometric structure,
at least on complete varieties. Our leitfaden is the question whether equally interesting
geometries exist in higher dimension.

We study here examples of rank $2$ vector bundles on \emph{abelian} threefolds,
partly for the intrinsic interest, and partly in the hope that the abelian
case may shed light on the case of general Calabi-Yau threefolds, but be
more accessible. Our central tool, besides the Serre construction, is monads:
these are usually put to work on rational varieties, and it
may be slightly surprising that they can be useful also in our context. On
the other hand, we do not know whether the bundles we construct, and their
moduli, show typical or exceptional behaviour.

\subsection{Notation}

We work over an algebraically closed field $k$ of characteristic zero. Our
terminology regarding (semi-) stable sheaves and their moduli follows Simpson
\cite{simpson94}; in particular, stability for a coherent $\OO_X$-module on a
polarized projective variety $(X,H)$ is defined using the normalized Hilbert
polynomial with respect to the polarization $H$. Stable sheaves admit a coarse
moduli space $M$, with a compactification $\overline{M}$ parametrizing
S-equivalence classes of semistable sheaves.
The sheaves we construct in this text will in fact have the stronger
property of $\mu$-stability, in the sense of Mumford and Takemoto,
which is measured by the slope, i.e.~the ratio of degree with respect to $H$
to the rank, and which implies stability. Conversely, semistability implies
$\mu$-semistability.

The words line bundles and vector bundles are used as synonyms of invertible
and locally free sheaves. In particular, an inclusion of vector bundles means
an inclusion as sheaves, i.e.~the quotient need not be locally free.
We take Chern classes to live in the Chow ring modulo numerical equivalence.

Let $(X,\Theta)$ be a principally polarized abelian variety. If $x\in X$
is a point, we write $T_x\colon X\to X$ for the translation map,
and define $\Theta_x$ as $T_x(\Theta) = \Theta+x$. We identify $X$ with its dual
$\Pic^0(X)$ by associating with $x$ the line bundle $\sheaf{P}_x =
\OO_X(\Theta-\Theta_x)$. The normalized Poincar\'e line bundle on
$X\times X$ is denoted $\sheaf{P}$; its restriction to $X\times\{x\}$ is
$\sheaf{P}_x$.

\section{The Serre construction}

In this section we apply the standard Serre construction to produce
rank $2$ vector bundles on principally polarized abelian threefolds,
including examples with small $c_2$. This is the content of Theorem
\ref{thm:existence}. These examples (in the case
of even $c_1$) will be our objects of study for the rest of this
paper.

\subsection{The bundles/curves correspondence}

Let $\sheaf{E}$ be rank $2$ vector bundle on a projective variety $X$,
and let $s\in\Gamma(X,\sheaf{E})$ be a section. If the vanishing locus
$V(s)$ has codimension $2$, then: (1) it is a locally complete intersection,
and (2) its canonical bundle is $\res{(\omega_X\tensor\Alt^2\sheaf{E})}{V(s)}$.
The Serre construction says (under a cohomological condition on $\Alt^2\sheaf{E}$)
that any codimension two subscheme $Y\subset X$
with these two properties is of the form $V(s)$. More precisely:

\begin{theorem}\label{thm:serre}
Let $X$ be a projective variety with a line bundle $\sheaf{L}$ satisfying
$H^p(X,\sheaf{L}^{-1})=0$ for $p=1,2$. Let $Y\subset X$ be a codimension two locally
complete intersection subscheme with canonical bundle isomorphic to
$\res{(\omega_X\tensor\sheaf{L})}{Y}$. Then there is a canonical isomorphism
\begin{equation*}
\Hom(\res{(\omega_X\tensor\sheaf{L})}{Y},\omega_Y) \iso \Ext^1(\sheaf{I}_Y\tensor \sheaf{L}, \OO_X)
\end{equation*}
which is functorial in $Y$ with respect to inclusions, and such that
isomorphisms on the left correspond to locally free
extensions on the right.
\end{theorem}

For the proof we refer to Hartshorne \cite[Thm.~1.1 and Rem.~1.1.1]{hartshorne78}, who attributes
``all essential ideas'' to Serre \cite{serre63}.

It follows that, whenever we choose an isomorphism
$\res{(\omega_X\tensor\sheaf{L})}{Y}\iso \omega_Y$,
the theorem gives an extension
\begin{equation}\label{eq:koszul}
\begin{diagram}
0 & \rTo & \OO_X & \rTo^{s} & \sheaf{E} & \rTo & \sheaf{I}_Y\tensor\sheaf{L} & \rTo & 0
\end{diagram}
\end{equation}
with $\sheaf{E}$ locally free, and hence $Y = V(s)$ as required.

\begin{definition}
We say that $\sheaf{E}$ and $Y$ \emph{corresponds} if there is a
short exact sequence \eqref{eq:koszul}.
\end{definition}

Note that, if $Y$ has several connected components, there may be several
non-isomorphic bundles $\sheaf{E}$ corresponding to $Y$. See Proposition
\ref{prop:serrefamily}.

\subsection{Construction of bundles}

For the rest of this paper, we fix a principally polarized abelian threefold
$(X,\Theta)$. We assume that its Picard number is one, although this assumption
is not essential in later sections. Thus
every divisor is numerically equivalent to an
integral multiple of $\Theta$. Moreover (see e.g.~Debarre \cite{debarre94}), an
application of the endomorphism construction of Morikawa \cite{morikawa54} and
Matsusaka \cite{matsusaka59} shows that every $1$-cycle is numerically
equivalent to an integral multiple of $\Theta^2/2$.

So fix classes $c_1 = m\Theta$ and $c_2=n\Theta^2/2$, where $m$ and $n$ are integers.
If these are the Chern classes of a rank two
vector bundle $\sheaf{E}$, then, by Riemann-Roch
\begin{equation*}
\chi(\sheaf{E}) = \tfrac{1}{6}(c_1^3-3c_1c_2) = m^3-\tfrac{3}{2}nm,
\end{equation*}
so either $m$ or $n$ is even. Moreover, if $\sheaf{E}$
is $\mu$-semistable, then Bogomolov's inequality reads $m^2 \le 2n$.

\begin{theorem}\label{thm:existence}
Let $(X, \Theta)$ be a principally polarized abelian threefold of Picard
number $1$, and let $c_1 = m\Theta$ and $c_2 = n\Theta^2/2$, with $m$ and $n$ integers.
Assume
\begin{enumerate}
\item the strict Bogomolov inequality holds, i.e.~$m^2<2n$, and
\item $n$ is even and $mn$ is divisible by $4$.
\end{enumerate}
Then there exist $\mu$-stable rank $2$ vector bundles with Chern classes
$c_1$ and $c_2$.
\end{theorem}

\begin{remark}
For each $c_1\in \NS(X)$, the theorem realizes every second $c_2$ that is allowed
by (strict) Bogomolov and Riemann-Roch. The other half seems much
more subtle. In fact, we do not know any example of a rank $2$ vector
bundle, stable or not, that violates condition (2). The situation
in which equality occurs in the Bogomolov inequality will be analysed
in Proposition \ref{prop:bogomolov-eq}.
\end{remark}

Before proving the theorem, we rephrase $\mu$-stability for $\sheaf{E}$
as a condition on the corresponding curve $Y$. The argument
is similar to that of Hartshorne \cite[Prop.~3.1]{hartshorne78} in
the case of $\PP^3$.

\begin{lemma}\label{lemma:stab}
Let $(X, \Theta)$ be as in the theorem, and $\sheaf{E}$ be a rank $2$
vector bundle
corresponding to a curve $Y\subset X$.
Let $c_1(\sheaf{E}) = m\Theta$. Then the following
are equivalent.
\begin{enumerate}
\item $\sheaf{E}$ is $\mu$-stable.
\item $m>0$ and $Y$ is not contained in any translate of any divisor in the linear
system $|k\Theta|$, where $k$ is the round down of $m/2$.
\end{enumerate}
\end{lemma}

\begin{proof}
Since $\sheaf{E}$ has a section, it is clear that $m>0$ is necessary for its
$\mu$-stability. Write $\floor{m/2}$ and $\ceil{m/2}$ for the round down and round up
of $m/2$. The bundle $\sheaf{E}$ fails $\mu$-stability if and only if it
contains a line bundle $\sheaf{P}_x(l\Theta)\subset\sheaf{E}$ with $l \ge m/2$. Since $\sheaf{P}_x(l\Theta)$ has global sections for $l$
positive, it suffices to test with $l=\ceil{m/2}$. Thus $\sheaf{E}$ is
$\mu$-stable if and only if
\begin{equation}\label{eq:vanish}
H^0(X, \sheaf{E}(-\ceil{m/2}\Theta)\tensor \sheaf{P}_x))=0
\quad\text{for all $x\in X$.}
\end{equation}
Now twist the short exact sequence \eqref{eq:koszul} with
$-\ceil{m/2}\Theta$ and take cohomology. Since $H^i(X,\OO_X(-\ceil{m/2}\Theta)))=0$ for
$i=0,1$, and the determinant of $\sheaf{E}$ has the form $\sheaf{P}_a(m\Theta)$ for
some $a\in X$, we find that the vanishing \eqref{eq:vanish} is equivalent to the
vanishing of $H^0(X, \sheaf{I}_Y(\floor{m/2}\Theta)\tensor \sheaf{P}_x)$ for all
$x\in X$. Since $\Theta$ is ample, this is equivalent to
\begin{equation*}
H^0(X, \sheaf{I}_Y\tensor T_x^*\OO_X(\floor{m/2}\Theta))=0
\quad\text{for all $x\in X$}
\end{equation*}
which is condition (2).
\end{proof}

\begin{proof}[Proof of Theorem \ref{thm:existence}]
Since $\mu$-stability, and the conditions (1) and (2) in the statement
of the theorem, are preserved under tensor product with line
bundles, it suffices to prove the theorem for $m=2$ and $m=3$.

When $m=2$, the theorem claims that there are $\mu$-stable rank $2$
bundles with $c_1 = 2\Theta$ and $c_2=N\Theta^2$ for all integers $N\ge2$. For
this, choose $N$ generic points $a_i\in X$ and let
\begin{equation*}
Y = \Union_{i=1}^N Y_i,\quad Y_i = \Theta_{a_i}\intersect\Theta_{-a_i}.
\end{equation*}
We want to apply the Serre construction to this curve.

First we claim that the $Y_i$'s are pairwise disjoint, for $a_i$ chosen
generically. In fact, for $i\ne j$ write
\begin{equation*}
Y_i\intersect Y_j =
\underbrace{(\Theta_{a_i}\intersect\Theta_{a_j})}_V \intersect
\underbrace{(\Theta_{-a_i}\intersect\Theta_{-a_j})}_W,
\end{equation*}
where $V$ and $W$ have codimension $2$.
By an easy moving lemma for abelian varieties \cite[Lemma 5.4.1]{BL92}, a general
translate $V+x$ intersects $W$ properly, hence empty. Thus (replacing $x$ by
a ``square root'' $x/2$) also $V+x$ and $W-x$ are disjoint. So $Y_i$ and
$Y_j$ will be disjoint after a small perturbation $a_i\mapsto a_i+x$,
$a_j\mapsto a_j+x$.

The normal bundle of each $Y_i\subset X$ is
$\OO_{Y_i}(\Theta_{a_i})\oplus\OO_{Y_i}(\Theta_{-a_i})$, hence the canonical
bundle $\omega_{Y_i}$ is $\OO_{Y_i}(\Theta_{a_i}+\Theta_{-a_i})$. The
theorem of the square shows that
$\Theta_{a_i}+\Theta_{-a_i}$ is linearly
equivalent to $2\Theta$. Since the $Y_i$'s are disjoint,
we conclude that $Y$ is a locally complete
intersection with canonical bundle $\OO_Y(2\Theta)$. The Serre
construction produces a bundle $\sheaf{E}$ with determinant
$\OO_X(2\Theta)$ and second Chern class $[Y] = \sum_i [Y_i] = N\Theta^2$.

Next we show $\mu$-stability. We claim that the only
theta-translates containing $Y_i$ are $\Theta_{a_i}$ and
$\Theta_{-a_i}$.
This is a standard result: the intersection of two theta-translates are never
contained in a third one. In fact, consider the Koszul complex:
\begin{equation*}
0 \to \OO_X(-\Theta_{a_i}-\Theta_{-a_i})
\to \OO_X(-\Theta_{a_i})\oplus\OO_X(-\Theta_{-a_i}) \to \sheaf{I}_{Y_i} \to 0.
\end{equation*}
Twist with an arbitrary theta-translate $\Theta_x$ and apply cohomology
to obtain an isomorphism
\begin{equation*}
H^0(X, \OO_X(\Theta_x - \Theta_{a_i}))
\oplus
H^0(X, \OO_X(\Theta_x - \Theta_{-a_i}))
\iso
H^0(X, \sheaf{I}_{Y_i}(\Theta_x)).
\end{equation*}
Thus $\Theta_x$ contains $Y_i$ if and only if $x=\pm a_i$ as claimed. It
follows that, for $N\ge 2$, no theta-translate contains $Y$, and so $\sheaf{E}$ is
$\mu$-stable by Lemma \ref{lemma:stab}.

In the case $m=3$, we take
\begin{equation*}
Y = \Union_{i=1}^N Y_i,\quad Y_i = D_i\intersect \Theta_{-2a_i}
\end{equation*}
for $N$ generic points $a_i\in X$ and generic divisors $D_i \in
|2\Theta_{a_i}|$. A similar argument to the one above shows that the Serre construction
produces a $\mu$-stable rank $2$ vector bundle with determinant $\OO_X(3\Theta)$ and
second Chern class $2N\Theta^2$, for each $N\ge 2$.
\end{proof}

Recall that a vector bundle $\sheaf{E}$ is \emph{semihomogeneous} if, for every
$x\in X$, there exists a line bundle $\sheaf{L}\in \Pic^0(X)$ such that
$T_x^*(\sheaf{E})$ is isomorphic to $\sheaf{E}\tensor\sheaf{L}$ (in short,
homogeneous means translation invariant, and semihomogeneous means translation
invariant up to twist). Semihomogeneous bundles are well understood thanks to
work of Mukai \cite{mukai78}.

\begin{proposition}\label{prop:bogomolov-eq}
Let $(X,\Theta)$ be as in the theorem, and let $c_1=m\Theta$ and $c_2=n\Theta^2/2$
satisfy $m^2 = 2n$, i.e.~equality occurs in the Bogomolov inequality.
Then $\sheaf{E}$ is a non simple, semihomogeneous vector bundle. In
particular, it is semistable, but not stable.

More precisely, there are a line bundle $\sheaf{L}_0$ and
points $x,y\in X$ such that $\sheaf{E}_0 = \sheaf{E}\tensor\sheaf{L}_0^{-1}$
is an extension (necessarily split if $x\ne y$)
\begin{equation*}
0 \to \sheaf{P}_x \to \sheaf{E}_0 \to \sheaf{P}_y \to 0.
\end{equation*}
\end{proposition}

\begin{proof}
Semihomogenous
bundles of rank $r$ are numerically characterized (Yang \cite{yang89}) by the property that the
Chern roots may be taken to be $c_1/r$. This means that the Chern character
takes the form $ch = r\exp(c_1/r)$, or, equivalently, the total Chern class is
$c = (1+c_1/r)^r$. If $r=2$, this is equivalent to $c_1^2=4c_2$. Thus
$\sheaf{E}$ is semihomogeneous.

Now we use several results by Mukai \cite{mukai78} on semihomogeneous bundles.
\emph{Simple} semihomogeneous vector bundles
are classified, up to twist by homogeneous line bundles, by the element $\delta = c_1/r$
in $\NS(X)\tensor\QQ$. But $m$ is even, since $m^2=2n$, so there exist line bundles
with class $c_1/2$, which rules out the possibility that $\sheaf{E}$ is simple.
Moreover, any semihomogeneous bundle is Gieseker-semistable, and it is simple
if and only if it is Gieseker-stable. This proves the first part.

For the last part, we use Mukai's Harder-Narasimhan filtration for semihomogeneous
bundles,
which in particular says that any semihomogeneous vector bundle with
$\delta=c_1/r$ has a filtration whose factors are simple semihomogeneous
bundles with the same invariant $\delta$. Choosing $\sheaf{L}_0$ in the
(integral) class $c_1/2$, we ensure that $\sheaf{E}_0$ has $\delta=0$.
Since $\sheaf{E}_0$ is semihomogeneous, but not simple, its
Harder-Narasimhan factors are necessarily line bundles with $c_1=0$.
\end{proof}

\subsection{A note on the curves $\Theta_a\intersect\Theta_{-a}$}

For later use, we make an observation regarding the curve obtained
by intersecting two general
theta-translates, which was used as input for the
Serre construction above (in the even $c_1$ case).

First note that there is a Zariski open subset $U\subset X$
such that $\Theta\intersect\Theta_x$ is a nonsingular
irreducible curve for all $x\in U$. This is standard: since $\Theta$
is a nonsingular surface, generic smoothness shows that
$\Theta\intersect\Theta_x$ is nonsingular, but possibly disconnected,
for generic $x$ (see Hartshorne \cite[III~10.8]{hartshorne77}). On
the other hand, $\Theta\intersect\Theta_x$ is an ample divisor on
$\Theta$, hence it is connected (see Hartshorne \cite[III~7.9]{hartshorne77}).

\begin{lemma}\label{lemma:theta-intersect}
Let $a$ and $b$ be two points in $X$ and define
$Y_a = \Theta_a\intersect\Theta_{-a}$ and $Y_b =
\Theta_b\intersect\Theta_{-b}$. Then, for $a$ and $b$ generic, 
no divisor in $|2\Theta|$ contains both $Y_a$ and $Y_b$.
\end{lemma}

\begin{proof}
Begin by imposing the conditions on $a$ and $b$ that $Y_a$ and $Y_b$ are
disjoint irreducible curves, and also that the two curves $\Theta_a\intersect\Theta_{\pm b}$ are
irreducible. Assume there is a divisor $D\in |2\Theta|$
containing both $Y_a$ and $Y_b$. We will prove the lemma by producing a
curve $C$ such that $C\intersect\Theta_b = C\intersect\Theta_{-b}$, and then
deduce from this that $b$ is not generic.

First we observe that $D$ meets $\Theta_a\intersect\Theta_b$ properly. As the
latter is irreducible, it suffices to verify that it is not contained in
$D$. In fact,
one checks (determine
$H^0(\sheaf{I}_{\Theta_a\intersect\Theta_b}(2\Theta))$ using the Koszul resolution) that the linear subsystem of
$|2\Theta|$, consisting of divisors containing $\Theta_a\intersect\Theta_b$, is
the pencil spanned by $\Theta_a+\Theta_{-a}$ and $\Theta_b+\Theta_{-b}$. The
only element of this pencil containing $Y_a$ is $\Theta_a+\Theta_{-a}$, and the
only element containing $Y_b$ is $\Theta_b+\Theta_{-b}$, so no element contains
both.

In particular, $D$ and $\Theta_a$ intersects properly, so 
$D\intersect\Theta_a$ is a curve containing $Y_a$. Since
$D\intersect\Theta_a$ has cohomology class $2\Theta^2$, and $Y_a$ has class
$\Theta^2$, there is another effective $1$-cycle $C$ of class $\Theta^2$
such that
\begin{equation*}
D\intersect\Theta_a = Y_a + C
\end{equation*}
as $1$-cycles. We saw above that $D\intersect\Theta_a$ meets $\Theta_b$ properly,
so we consider the $0$-cycle
\begin{equation*}
D\intersect\Theta_a\intersect\Theta_b = Y_a\intersect\Theta_b + C\intersect\Theta_b.
\end{equation*}
The left hand side contains $Y_b\intersect\Theta_a$. Since $Y_a$ and $Y_b$ are
disjoint, this means that $C\intersect\Theta_b$ contains $Y_b\intersect\Theta_a$,
i.e.~their difference is an effective cycle. But these are $0$-cycles
of the same degree, so they are equal. None of the arguments given distinguish
between $b$ and $-b$, so we find that also $C\intersect\Theta_{-b}$ equals
$Y_b\intersect\Theta_a$. Thus we have established
\begin{equation*}
C\intersect\Theta_{b} = C\intersect\Theta_{-b}.
\end{equation*}

To conclude, we apply the endomorphism construction of Morikawa \cite{morikawa54}
and Matsusaka \cite{matsusaka59}, which we briefly recall.
The endomorphism $\alpha = \alpha(C,\Theta)$ associated to $C$ and $\Theta$ is
defined by
\begin{equation*}
\alpha(x) = \sum(C\cdot\Theta_x) - \sum(C\cdot\Theta)
\end{equation*}
where each term means the sum, using the group law, of the points
in the intersection cycle appearing. This is well 
defined as a point in $X$,
although the intersection cycle is only defined up to rational equivalence.
The constant term is included to force $\alpha(0)=0$, i.e.~to make $\alpha$ a
group homomorphism. We have just established that $C$ intersects $\Theta_{b}$
and $\Theta_{-b}$ properly, and the two intersections are equal already as cycles.
In particular $\alpha(b)=\alpha(-b)$, so all we need to know to prove the lemma
is that $\alpha$ is not constant, so that $\alpha(2b) \ne 0$ defines a nonempty
Zariski open subset. But in fact, a theorem of
Matsusaka \cite{matsusaka59} tells us
that $\alpha$ is multiplication by $2$ (the intersection number $C\cdot\Theta=3!$
divided by $\dim X=3$), so the condition required is just that $4b\ne 0$, i.e.~$b$
is not a $4$-torsion point.
\end{proof}

As an immediate consequence of the Lemma, we find that if $\sheaf{E}$
corresponds to a curve with at least two components of the form $\Theta_{a_i}\intersect\Theta_{-a_i}$,
for sufficiently general points $a_i$, then the short exact sequence
\begin{equation*}
\begin{diagram}
0 & \rTo & \OO_X & \rTo^{s} & \sheaf{E} & \rTo & \sheaf{I}_Y(2\Theta) & \rTo & 0,
\end{diagram}
\end{equation*}
shows that
$H^0(X,\sheaf{E})$ is spanned by $s$.

\begin{proposition}\label{prop:serrefamily}
Fix $N\ge 2$ general points $a_i\in X$, and let $Y$ be the union
of the curves $Y_i = \Theta_{a_i}\intersect\Theta_{-a_i}$. Then
the vector bundles
$\sheaf{E}$ corresponding to the union $Y$ of
$Y_i=\Theta_{a_i}\intersect\Theta_{-a_i}$ form an $(N-1)$-dimensional family,
parametrized by $\Gm^{N-1}$.
\end{proposition}

\begin{proof}
The Serre construction gives a one-one correspondence between isomorphisms
$\omega_Y\iso\OO_Y(2\Theta)$ modulo scale, and isomorphism classes of
vector bundles $\sheaf{E}$ which admit a section vanishing at $Y$: the choice
of a section can be left out, since we just observed that it is unique modulo
scale. But isomorphisms $\omega_Y\iso\OO_Y(2\Theta)$ constitute a homogeneous
$\Gm^N$-space, as $Y$ has $N$ connected components. Dividing by scale, we
are left with $\Gm^N/\Gm\iso \Gm^{N-1}$.
\end{proof}

\section{Monads}

It turns out that the vector bundles constructed in Theorem \ref{thm:existence}
admit more deformations than are visible in the Serre construction, i.e.~more
deformations than those obtained by varying the curve $Y$ and the isomorphism
$\omega_Y\iso\OO_X(2\Theta)$. In this section we rephrase the construction in
terms of certain monads (Proposition \ref{prop:decomposable}).  This new viewpoint is
then used in the remaining sections to analyse first order deformations.

\begin{definition}[Barth--Hulek \cite{BH78}]
A \emph{monad} is a composable pair of maps of vector bundles
\begin{equation*}
\begin{diagram}
\sheaf{A} & \rTo^{\phi} & \sheaf{B} & \rTo^{\psi} & \sheaf{C}
\end{diagram}
\end{equation*}
such that $\psi\circ\phi$ is zero, $\psi$ is surjective and $\phi$ is an embedding of
vector bundles (i.e.~injective as a homomorphism of sheaves, and
with locally free cokernel).

Thus $\sheaf{E}=\Ker(\psi)/\Image(\phi)$ is a vector bundle, and we say
that the monad is a \emph{monad for $\sheaf{E}$}.
\end{definition}

We will also use chain complex notation $(M^\spot, d)$ for monads,
so that $M^{-1}=\sheaf{A}$, $M^0=\sheaf{B}$, $M^1=\sheaf{C}$ and $M^i$ is zero otherwise, and
the differential $d$ consists of two nonzero components $d^{-1}=\phi$ and
$d^0=\psi$. Thus $M^\spot$ is exact except in degree zero, where its
cohomology is $\sheaf{E}=H^0(M^\spot)$.

\subsection{Decomposable monads}

Consider rank $2$ vector bundles $\sheaf{E}$ with trivial determinant $\Alt^2\sheaf{E}\iso\OO_X$
on the principally polarized abelian threefold $(X,\Theta)$. 
From the construction in
Theorem \ref{thm:existence}, we have a series of such vector bundles, such that
$\sheaf{E}(\Theta)$ corresponds to a curve $Y=\Union_i Y_i$, where
$Y_i=\Theta_{a_i}\intersect\Theta_{-a_i}$. (The assumption that $X$ has
Picard number $1$ is not needed here; this was only needed to establish $\mu$-stability of $\sheaf{E}$,
which is not relevant in this section.)

We now show that, corresponding to the decomposition of $Y$ into its
connected components $Y_i$, there is a way of building up $\sheaf{E}$
from the Koszul complexes\footnote{Here and elsewhere, whenever $f\colon \sheaf{F}_1\to\sheaf{F}_2$
is a homomorphism of sheaves, we use the same symbol to denote any twist $f\colon \sheaf{F}_1(D)\to\sheaf{F}_2(D)$.}
\begin{equation}\label{eq:koszul-i}
\begin{diagram}
\xi_i\colon&
0 & \rTo & \OO_X(-\Theta) & \rTo^\koszA & \sheaf{P}_{a_i}\plus \sheaf{P}_{-a_i} & \rTo^\koszB & \sheaf{I}_{Y_i}(\Theta) & \rTo & 0
\end{diagram}
\end{equation}
where $\vartheta^\pm_i$ are nonzero global sections of $\OO_X(\Theta_{\pm a_i})$.
This can be conveniently phrased in terms of a monad.

\begin{proposition}\label{prop:decomposable}
Let $a_1,\dots,a_N\in X$ be generically chosen points and
$Y_i=\Theta_{a_i}\intersect\Theta_{-a_i}$. Then $\sheaf{E}(\Theta)$ corresponds
to $Y=\Union_{i=1}^N Y_i$ if and only if $\sheaf{E}$ is isomorphic to the cohomology of
a monad
\begin{equation*}
\begin{diagram}
(N-1)\OO_X(-\Theta) &
\rTo^{\phi} &
\Plus_{i=1}^N (\sheaf{P}_{a_i}\plus \sheaf{P}_{-a_i}) &
\rTo^{\psi} &
(N-1)\OO_X(\Theta)
\end{diagram}
\end{equation*}
where, if we decompose $\phi$ and $\psi$ into pairs
\begin{align*}
\phi^\pm\colon& (N-1)\OO_X(-\Theta) \to \Plus_{i=1}^N \sheaf{P}_{\pm a_i}\\
\psi^\pm\colon& \Plus_{i=1}^N \sheaf{P}_{\pm a_i} \to (N-1)\OO_X(\Theta)
\end{align*}
then we have
\begin{equation*}
\phi^{\pm} =
\left(\begin{matrix}
\vartheta_1^\pm \\
& \vartheta_2^\pm \\
& & \ddots\\
& & & \vartheta_{N-1}^\pm\\
\vartheta_N^\pm & \vartheta_N^\pm & \cdots & \vartheta_N^\pm
\end{matrix}\right),
\quad
\psi^\pm = \pm \dual{(\phi^\mp)}
\end{equation*}
for nonzero sections $\vartheta^\pm_i\in\Gamma(X,\OO_X(\Theta_{\pm a_i}))$.
\end{proposition}

\begin{proof}
One immediately verifies that homomorphisms $\phi$ and $\psi$ of this form do define
a monad.

The statement that $\sheaf{E}(\Theta)$ and $Y$ correspond means that $\sheaf{E}$
is an extension
\begin{equation*}
\begin{diagram}
\xi\colon& 0 & \rTo & \OO_X(-\Theta) & \rTo & \sheaf{E} & \rTo & \sheaf{I}_Y(\Theta) & \rTo & 0.
\end{diagram}
\end{equation*}
Giving such an extension is, by Theorem \ref{thm:serre}, equivalent to
giving an isomorphism $\OO_Y(2\Theta)\iso\omega_Y$. The obvious
decomposition
\begin{equation*}
\Hom(\OO_Y(2\Theta), \omega_Y) \iso \Plus_{i=1}^N \Hom(\OO_{Y_i}(2\Theta), \omega_{Y_i})
\end{equation*}
gives, when applying Theorem \ref{thm:serre} also to each $Y_i$,
a corresponding decomposition
\begin{equation}\label{eq:ext-decomp}
\Ext^1(\sheaf{I}_Y(\Theta), \OO_X(-\Theta)) \iso \Plus_{i=1}^N \Ext^1(\sheaf{I}_{Y_i}(\Theta), \OO_X(-\Theta)),
\end{equation}
which sends $\xi$ to an $N$-tuple of extensions $\xi_i$. Each
$\Hom(\OO_{Y_i}(2\Theta), \omega_{Y_i})$ is one dimensional, since $Y_i$ is
connected, so $\Ext^1(\sheaf{I}_{Y_i}(\Theta), \OO_X(-\Theta))$ is one
dimensional, too. This shows that each $\xi_i$ is of the form
\eqref{eq:koszul-i}.

From the functoriality in Theorem \ref{thm:serre},
it follows that the inclusion of each
direct summand in \eqref{eq:ext-decomp} is the natural map, induced by the
inclusion $\sheaf{I}_Y\subset\sheaf{I}_{Y_i}$. Thus $\xi$ is obtained from the
$\xi_i$'s by pulling them back over this inclusion of ideals, and adding the
results in $\Ext^1(\sheaf{I}_Y(\Theta), \OO_X(-\Theta))$. By definition of
(Baer) addition in Ext-groups, this means that there is a commutative diagram
\begin{equation*}
\begin{diagram}
0 & \rTo & N\OO_X(-\Theta) & \rTo & \Plus_{i=1}^N (\sheaf{P}_{a_i}\plus \sheaf{P}_{-a_i}) & \rTo & \Plus_{i=1}^N\sheaf{I}_{Y_i}(\Theta) & \rTo & 0 \\
 && \dOnto^\beta && \dOnto && \dEq \\
0 & \rTo & \OO_X(-\Theta) & \rTo & \sheaf{F} & \rTo & \Plus_{i=1}^N\sheaf{I}_{Y_i}(\Theta) & \rTo & 0 \\
 && \uEq && \uInto && \uInto^{\alpha} \\
0 & \rTo & \OO_X(-\Theta) & \rTo & \sheaf{E} & \rTo & \sheaf{I}_Y(\Theta) & \rTo & 0
\end{diagram}
\end{equation*}
where the top row is $\Plus_i\xi_i$, the bottom row is $\xi$, the top left square is pushout over the
$N$-fold addition $\beta$, the bottom right square is pullback along the
inclusion $\alpha$, and $\sheaf{F}$ is just an intermediate sheaf (in fact a
vector bundle) that we do not care about. This diagram presents $\sheaf{E}$
as the middle cohomology of a complex
\begin{equation*}
\begin{diagram}
\Ker(\beta) & \rTo^\phi & \Plus_{i=1}^N (\sheaf{P}_{a_i}\plus \sheaf{P}_{-a_i}) & \rTo^\psi \Coker(\alpha).
\end{diagram}
\end{equation*}

Now identify $\Ker(\beta)$ with $(N-1)\OO_X(-\Theta)$ by means of the monomorphism
\begin{equation*}
(N-1)\OO_X \to N\OO_X,\quad
(f_1,\dots,f_{N-1})\mapsto (f_1,\dots,f_{N-1},-\tsum_i f_i)
\end{equation*}
and similarly identify $\Coker(\alpha)$ with $(N-1)\OO_X(\Theta)$ by means of
the epimorphism
\begin{equation*}
N\OO_X \to (N-1)\OO_X,\quad
(f_1,\dots,f_N) \mapsto (f_1-f_N,\dots,f_{N-1}-f_N)
\end{equation*}
(the latter is surjective even when restricted to
$\Plus_i\sheaf{I}_{Y_i}$ because the $Y_i$'s are pairwise disjoint).
Via these identifications, the homomorphisms $\phi$ and $\psi$ are represented
by the matrices as claimed, except that $\vartheta^\pm_N$ appears with opposite sign.
Change its sign, and we are done.
\end{proof}

\begin{definition}
A monad is \emph{decomposable} if it is isomorphic, as a complex, to
a monad of the form appearing in Proposition \ref{prop:decomposable}.
\end{definition}

With this terminology, a rank $2$ vector bundle $\sheaf{E}$ can be
resolved by a decomposable monad if and only if $\sheaf{E}(\Theta)$
corresponds to a disjoint union $Y=\Union_i Y_i$, where $Y_i =
\Theta_{a_i}\intersect\Theta_{-a_i}$, via the Serre construction.

\begin{remark}
The symmetry seen in the decomposable monads is no accident, but
reflects the self duality of $\sheaf{E}$ corresponding to the
natural pairing $\alt$ on $\sheaf{E}$ with values
in $\Alt^2(\sheaf{E})\iso\OO_X$. See Barth--Hulek \cite{BH78}.
\end{remark}

\section{Digression on the hyperext spectral sequence}

Our basic aim is to understand first order deformations of the bundles
$\sheaf{E}$ appearing as the cohomology of a decomposable monad.  The strategy
is to analyse $\Ext^1(\sheaf{E},\sheaf{E})$ using the first hyperext spectral
sequence associated to the monad. This is in principle straight forward, but
requires some honest calculation. As preparation, we collect in this section a
few standard constructions in homological algebra, for ease of reference.
We fix an abelian category $\cat{A}$ with enough injectives and infinite
direct sums, and denote by $K(\cat{A})$ the homotopy category of complexes
and by $D(\cat{A})$ the derived category.

\subsection{The spectral sequence}\label{sec:spectral}

Let $(M^\spot, d_M)$ and $(N^\spot, d_N)$ denote complexes
in $\cat{A}$, and assume that $N^\spot$ is bounded from
below. The \emph{first hyperext spectral sequence} is a spectral
sequence
\begin{equation*}
E^{pq}_1 = \Plus_i \Ext^q(M^i,N^{i+p}) \abuts \Ext^{p+q}(M^\spot, N^\spot).
\end{equation*}
Briefly, take a double injective resolution
$N^\spot \to I^{\spot\spot}$ with $I^{\spot\spot}$ concentrated
in the upper half plane (for instance a Cartan-Eilenberg resolution),
and form the double complex $\Hom^{\spot\spot}(M^\spot, I^{\spot\spot})$.
The required spectral sequence is the first spectral sequence associated
to this double complex.

\subsection{The edge map}\label{sec:edge}

Along the axis $q=0$, the first sheet of the spectral sequence in Section
\ref{sec:spectral} has the usual hom-complex $\Hom^\spot(M^\spot, N^\spot)$.
Its cohomology is
\begin{equation*}
E^{p,0}_2 = \Hom_{K(\cat{A})}(M^\spot, N^\spot[p])
\end{equation*}
where the right hand side denotes homotopy classes of morphisms
of complexes. Since all differentials
emanating from $E^{p,0}_r$ for $r\ge 2$ vanish, there are canonical
\emph{edge maps}
\begin{equation*}
E^{p,0}_2 \onto E^{p,0}_\infty \subset \Ext^p(M^\spot, N^\spot)
\end{equation*}
Viewing the right hand side as the group $\Hom_{D(\cat{A})}(M^\spot,
N^\spot[p])$ of morphisms in the derived category, it is reasonable to expect,
and not hard to verify, that the edge map is in fact the canonical
map
\begin{equation*}
\Hom_{K(\cat{A})}(M^\spot, N^\spot[p]) \to
\Hom_{D(\cat{A})}(M^\spot, N^\spot[p]).
\end{equation*}

Thus the image of $E^{p,0}_2$ in the limit object $\Ext^p(M^\spot, N^\spot)$,
consists of those $p$-extensions that can be
realized by actual morphisms $M^\spot\to N^\spot[p]$ between complexes, without
inverting quasi-isomorphisms.

\subsection{Differentials at $E_2$}
For $q=1$, it is convenient to view elements of $\Ext^1(M^i, N^{i+p})$
as extensions, in the sense of short exact sequences, and this
viewpoint leads to the following interpretation
of the differentials $d^{p1}_2$
at the $E_2$-level:

\begin{lemma}\label{lemma:obstruction}
Let $\xi \in E^{p1}_1$ be given as a collection of extensions
\begin{equation*}
\xi_i\colon 0 \to N^{i+p} \to X^i \to M^i \to 0.
\end{equation*}
\begin{enumerate}
\item We have $d^{p1}_1(\xi)=0$ if and only if there are maps $f^i$
such that the diagram
\begin{equation*}
\begin{diagram}
\cdots & \rTo & N^{i+p-1}         & \rTo^{d_N}         & N^{i+p}         & \rTo^{d_N}         & N^{i+p+1}         & \rTo & \cdots \\
       &      & \dTo        &                    & \dTo            &                        & \dTo \\
\cdots & \rTo & X^{i-1} & \rTo^{f^{i-1}} & X^i & \rTo^{f^i} & X^{i+1} & \rTo & \cdots \\
       &      & \dTo        &                    & \dTo            &                        & \dTo \\
\cdots & \rTo & M^{i-1}         & \rTo^{d_M}         & M^i         & \rTo^{d_M}         & M^{i+1}         & \rTo & \cdots \\
\end{diagram}
\end{equation*}
commutes.
\item If we have such a collection of maps $(f^i)$, then $\xi$ represents
an element of $E^{p1}_2$, and the differential
\begin{equation*}
d^{p1}_2\colon E^{p1}_2\to E^{p+2,0}_2 = \Hom_{K(\cat{A})}(M^\spot, N^\spot[p+2])
\end{equation*}
sends $\xi$ to the morphism having components $M^{i-1}\to N^{i+p+1}$ induced
by $f^i\circ f^{i-1}$. In particular $d^{p1}_2(\xi)=0$ if and only
if there exists a collection $(f^i)$ making the middle row in the diagram in (1)
a complex.
\end{enumerate}
\end{lemma}

\begin{proof}
This is straight forward, although tedious, to verify directly
from the construction of the spectral sequence.
\end{proof}

\subsection{Serre duality}\label{sec:duality}

Let $X$ be a scheme of pure dimension $d$ over a field, with a dualizing
sheaf $\omega_X$ such that Grothendieck-Serre duality holds.  Let $M^\spot$ be
a bounded below complex of coherent $\OO_X$-modules. We obtain two spectral
sequences from \eqref{eq:spectral}: one abutting to $\Ext^n(\OO_X, M^\spot) = H^n(X, M^\spot)$,
which we denote by $E$, and one
abutting to $\Ext^n(M^\spot, \omega_X)$, which we denote by $\hat{E}$.
Then $E$ is nothing but the first hypercohomology spectral sequence, and
the $E_1$-levels of $E$ and $\hat{E}$ are Grothendieck-Serre dual.
We need to know that the duality extends
to all sheets.

\begin{lemma}
The two spectral sequences $E$ and $\hat{E}$ are
dual in the following sense:
\begin{enumerate}
\item There are canonical dualities between the vector spaces
$E^{pq}_r$ and $\hat{E}^{-p, d-q}_r$ for all $p,q,r$, extending the
Grothendieck-Serre duality between $H^q(X,M^p)$
and $\Ext^{n-q}(M^p, \omega)$ for $r=1$.
\item The differentials
\begin{align*}
d^{pq}_r&\colon E^{pq}_r \to E^{p+r,q-r+1}_r\\
\hat{d}^{-p-r,d-q+r-1}_r&\colon \hat{E}^{-p-r,d-q+r-1}_r \to \hat{E}^{-p,d-q}_r
\end{align*}
are dual maps.
\end{enumerate}
\end{lemma}

(The statement can be extended to give a full duality
between the two spectral sequences, including the filtrations
on the abutments and all maps involved. The above
is sufficient for our needs.)

\begin{proof}
This seems to be well known. We include a sketch, following Herrera--Liebermann
\cite{HL71} (they work in a context where the complexes have differentials that are
differential operators of degree one; this demands more care than in our
situation). Firstly, for any
three complexes $L^\spot$, $M^\spot$, $N^\spot$, the Yoneda pairing
\begin{equation*}
\Ext^i(L^\spot, M^\spot) \times \Ext^j(M^\spot, N^\spot) \to \Ext^{i+j}(L^\spot, N^\spot)
\end{equation*}
can be defined on hyperext groups by resolving $M^\spot$ and $N^\spot$
by injective double complexes, and taking the double hom complex. On this
``resolved'' level, the Yoneda pairing is given by composition, and there is an
induced pairing of hyperext spectral sequences in the appropriate sense, which
specializes to the usual Yoneda pairing between ext groups of the individual
objects $L^l$, $M^m$, $N^n$ at the $E_1$-level. Specialize to the situation
$L^\spot=\OO_X$ and $N^\spot=\omega_X$ to obtain a morphism of spectral
sequences from $E$ to the dual of $\hat{E}$, in the above sense.
At the $E_1$-level this is the
Grothendieck-Serre duality map, hence an isomorphism, which is enough to conclude
that it is an isomorphism of spectral sequences \cite[Section 11.1.2]{EGAIII-1}.
\end{proof}

\begin{remark}
If $M^\spot$ and $N^\spot$ denote two complexes of vector bundles,
then we may apply the Lemma to the complex
$\dual{(M^\spot)} \tensor N^\spot$ to obtain a duality
between the two hyperext spectral sequences abutting to $\Ext^n(M^\spot,N^\spot)$
and $\Ext^n(N^\spot, M^\spot\tensor\omega_X)$, respectively.
\end{remark}

\section{Deformations of decomposable monads}

We now apply the homological algebra from the previous section to analyse
first order deformations of vector bundles $\sheaf{E}$ which can be resolved
by a decomposable monad. Firstly, we find that deformations obtained by varying
the isomorphism $\omega_Y\iso \OO_Y(2\Theta)$ in the Serre construction
coincides with the deformations obtained by varying the differential in
the monad, while keeping the objects fixed. Secondly, and this is the nontrivial
part, we find that all first order deformations of $\sheaf{E}$ can be obtained
by also deforming the objects in the monad, and there are more of these
deformations than those obtained by varying $Y$ in the Serre construction.
Since the objects in the monad are sums of line bundles, their first order
deformations are easy to understand, so we are able to compute the
dimension of $\Ext^1(\sheaf{E},\sheaf{E})$, in Theorem \ref{thm:quantitative}.

\subsection{Calculations in the spectral sequence}

Let $\sheaf{E}$ be the rank $2$ vector bundle given
as the cohomology of a decomposable monad
\begin{equation*}
\begin{diagram}
M^\spot\colon & \sheaf{A} & \rTo^\phi & \sheaf{B} & \rTo^{\psi} & \sheaf{C}
\end{diagram}
\end{equation*}
given explicitly in Proposition \ref{prop:decomposable}. In particular,
$\sheaf{C}=\dual{\sheaf{A}}$, and $\sheaf{B}$ is self dual. If we fix
the self duality $\iota\colon \sheaf{B}\to\dual{\sheaf{B}}$, given by the direct sum of the skew symmetric
\begin{equation*}
\begin{pmatrix}
0 & -1\\
1 & 0
\end{pmatrix}\colon
\sheaf{P}_{a_i}\plus\sheaf{P}_{-a_i} \to \sheaf{P}_{-a_i}\plus\sheaf{P}_{a_i},
\end{equation*}
then $\psi=\dual{\phi}\circ\iota$. More generally,
for any map $f\colon \sheaf{A}\to\sheaf{B}$, we define its transpose $f^t\colon \sheaf{B}\to\sheaf{C}$ by
\begin{equation*}
f^t = \dual{f}\circ\iota.
\end{equation*}
Thus $\psi$ is the transpose of $\phi$.\footnote{One can show
that, in the affine space of all homomorphisms $f\colon\sheaf{A}\to\sheaf{B}$,
the locally closed subset $U$ defined by (1) $f$ is an embedding of vector
bundles, and (2) the composition $f^t\circ f$ is zero, has an irreducible
connected component corresponding to decomposable monads. It seems plausible
that this component is all of $U$.}

The spectral sequence from Section \ref{sec:spectral} gives
\begin{equation}\label{eq:spectral}
E^{pq}_1 = \Plus_i \Ext^q(M^i, M^{i+p}) \abuts
\Ext^{p+q}(\sheaf{E},\sheaf{E}).
\end{equation}
Using that, for any $x\in X$, line bundles of the form $\sheaf{P}_x(m\Theta)$
have sheaf cohomology concentrated in degree $0$ when $m>0$ and in top
degree when $m<0$, we see that the nonzero terms in the first sheet
have the shape depicted in Figure \ref{fig:spectral}.
It follows that all differentials at level $E_r$
vanish for $r=3$ and $r>4$. Also, the duality of Section \ref{sec:duality},
applied to $M^\spot\tensor\dual{(M^\spot)}$, shows that each term $E^{pq}_r$
is dual to $E^{-p, 3-q}_r$, and similarly for the differentials.
In this section we analyse the $E_2$-sheet,
and get as a consequence that the spectral sequence in fact degenerates
at the $E_3$-level.

\begin{figure}
\begin{diagram}[height=2em]
E^{-2,3}_1 & \rTo & E^{-1,3}_1 & \rTo & E^{0,3}_1 \\
           &      &            &      & E^{0,2}_1 \\
           &      &            &      & E^{0,1}_1 \\
           &      &            &      & E^{0,0}_1 & \rTo & E^{1,0}_1 & \rTo & E^{2,0}_1
\end{diagram}
\caption{The first sheet in the spectral sequence for
$\Ext^i(\sheaf{E},\sheaf{E})$}\label{fig:spectral}
\end{figure}

\subsubsection{The objects $E^{pq}_2$}
By duality, it suffices to consider the lower
half of Figure \ref{fig:spectral}. The only nonzero differentials in
this area, at the $E_1$-level, are in the lower row $q=0$. We observed
in Section \ref{sec:edge} that the cohomology groups of this row
are the groups of morphisms $M^\spot\to M^\spot[p]$ modulo homotopy.

\begin{lemma}\label{lem:htpy}
The dimensions of $E^{p,0}_2$
for $p=0,1,2$ are $1$, $N-1$ and $6(N-1)^2-N+2$, respectively.
\end{lemma}

\begin{proof}
The vector spaces in question are the cohomologies of the complex
\begin{equation*}
\begin{diagram}
0 &\rTo & E^{0,0}_1 & \rTo^{d^{0,0}_1} & E^{1,0}_1 & \rTo^{d^{1,0}_1} &
E^{2,0}_1 & \rTo & 0,
\end{diagram}
\end{equation*}
where
\begin{equation}\label{eq:E1-dim}
\begin{split}
\dim E^{0,0}_1 &= \dim \big(
\Hom(\sheaf{A},\sheaf{A})\plus\Hom(\sheaf{B},\sheaf{B})\plus\Hom(\sheaf{C},\sheaf{C})\big)\\
&= 2(N-1)^2 + 2N\\
\dim E^{1,0}_1 &= \dim \big(\Hom(\sheaf{A}, \sheaf{B})\plus\Hom(\sheaf{B}, \sheaf{C})\big)\\
&= 4N(N-1)\\
\dim E^{2,0}_1 &= \dim \Hom(\sheaf{A}, \sheaf{C})\\
&= 8(N-1)^2
\end{split}
\end{equation}
(using that the space of global sections of $\OO_X(\Theta_{\pm a_i})$
has dimension $1$, and the space of global sections of $\OO_X(2\Theta)$
has dimension $8$).
Thus
it suffices to compute the dimensions of the kernels of the two
differentials $d^{0,0}_1$ and $d^{1,0}_1$, i.e.~the vector spaces of
morphisms of degree $0$ and $1$ from the monad to itself.

One checks immediately that any morphism
$M^\spot\to M^\spot$ (of degree $0$) is multiplication with a scalar,
so
\begin{equation}\label{eq:E002-dim}
\dim E^{0,0}_2 = 1.
\end{equation}

Next we compute the dimension of the space of morphisms $M^\spot\to
M^\spot[1]$. Since $\sheaf{C}=\dual{\sheaf{A}}$, such a morphism is given by an element of
\begin{equation*}
\Hom(\sheaf{A},\sheaf{B})\plus\Hom(\sheaf{B},\dual{\sheaf{A}}),
\end{equation*}
which we may write as $(\mu, -\nu^t)$, where both $\mu$ and $\nu$ are
homomorphisms $\sheaf{A}\to \sheaf{B}$. The sign on $-\nu^t$ is inserted
to compensate for the sign on the differential in the shifted complex
$M^\spot[1]$; thus $(\mu,-\nu^t)$ defines a morphism $M^\spot\to M^\spot[1]$
if and only if $\nu^t\circ\phi = \phi^t\circ\mu$.

As in Proposition \ref{prop:decomposable}, we
decompose these homomorphisms into pairs $\mu^\pm$ and $\nu^\pm$, and then
\begin{equation}\label{eq:compositions}
\begin{split}
\nu^t \circ \phi &= \dual{(\nu^-)}\circ\phi^+ - \dual{(\nu)^+}\circ\phi^- \\
\phi^t \circ \mu &= \dual{(\phi^-)}\circ\mu^+ -
\dual{(\phi^+)}\circ\mu^-.
\end{split}
\end{equation}
Choosing generators $\vartheta^\pm_i\in
\Gamma(X, \OO_X(\Theta_{\pm a_i}))$, we may represent $\mu$ by a matrix
with entries $\mu^\pm_{ij}\vartheta^\pm_i$, where $\mu^\pm_{ij}$ are
scalars. Similarly for $\nu$. Then the two compositions
\eqref{eq:compositions} are given by $(N-1)\times(N-1)$ scalar matrices with entries
\begin{equation}\label{eq:components}
\begin{split}
(\nu^t\circ\phi)_{ij} &= (\mu^+_{ij} - \mu^-_{ij})\vartheta^+_i\vartheta^-_i + (\mu^+_{Nj} -
\mu^-_{Nj})\vartheta^+_N\vartheta^-_N \\
(\phi^t\circ\mu)_{ij} &= (\nu^+_{ji} -
\nu^-_{ji})\vartheta^+_j\vartheta^-_j + (\nu^+_{Ni} -
\nu^-_{Ni})\vartheta^+_N\vartheta^-_N.
\end{split}
\end{equation}
Recall that the Kummer map $X \to |2\Theta|$ sends $a_i\in X$
to the divisor $\Theta_{a_i} + \Theta_{-a_i}$. This implies that,
for sufficiently general points $a_i$, and $i\ne j$, the three elements $\vartheta^+_i\vartheta^-_i$,
$\vartheta^+_j\vartheta^-_j$ and $\vartheta^+_N\vartheta^-_N$ are
linearly independent in $\Gamma(X, \OO_X(2\Theta))$. It follows easily that
the two expressions in \eqref{eq:components} coincide for all $i$ and
$j$ if and only if there are equalities of scalar $(N-1)\times (N-1)$
matrices
\begin{equation*}
(\mu^+_{ij}) - (\mu^-_{ij}) = (\nu^+_{ij}) - (\nu^-_{ij}) =
\left(\begin{matrix}
c_1\\
& c_2 \\
& & \ddots\\
& & & c_{N-1}\\
c_N & c_N & \cdots & c_N
\end{matrix}\right),
\end{equation*}
where $c_1,\dots, c_N$ are arbitrary scalars. Thus the vector space of
morphisms $M^\spot\to M^\spot[1]$ has a basis corresponding to the
$(\mu^+_{ij})$, $(\nu^+_{ij})$ and $(c_i)$, hence has dimension
$2N(N-1) + N$. The
expressions for $\dim E^{p,0}_2$ follow from this, together with
\eqref{eq:E1-dim} and \eqref{eq:E002-dim}.
\end{proof}

\subsubsection{The differentials $d^{pq}_2$}\label{sec:E2diff}
The only nonzero differentials at the $E_2$-level are $d^{0,1}_2$ and
its dual $d^{-2,3}_2$. So it suffices to analyse $d^{0,1}_2$. This is,
by Lemma \ref{lemma:obstruction}, an obstruction map for
equipping first order infinitesimal deformations of the objects $M^i$
with differentials, and will henceforth be denoted $\ob$.

The domain
\begin{equation}\label{eq:E012}
E^{0,1}_2 = \Plus_i \Ext^1(M^i, M^i)
\end{equation}
of $\ob= d^{0,1}_2$ is canonically isomorphic to a direct sum of a large number of
copies of $H^1(X, \OO_X)$. More precisely, for each $i$ and $j$ from $1$ to $N-1$,
apply the bifunctor $\Ext^1(-,-)$ to the $i$'th projection $\sheaf{A}\to \OO_X(-\Theta)$
in the first argument and the $j$'th inclusion $\OO_X(-\Theta)\to \sheaf{A}$ in the second
argument. This defines an inclusion
\begin{align*}
f_{ij}&\colon H^1(X,\OO_X)\iso \Ext^1(\OO_X(-\Theta),\OO_X(-\Theta))
\into \Ext^1(\sheaf{A}, \sheaf{A})
\intertext{and clearly the direct sum of all the $f_{ij}$'s is an isomorphism.
Similarly, for all $i$ and $j$ from $1$ to $N-1$, we define inclusions }
h_{ij}&\colon H^1(X,\OO_X)\iso \Ext^1(\OO_X(\Theta),\OO_X(\Theta))\into \Ext^1(\sheaf{C},\sheaf{C})
\intertext{whose direct sum is an isomorphism. Finally,
for all $i$ from $1$ to $N$, and each sign $\pm$, define inclusions}
g_i^\pm&\colon H^1(X,\OO_X)\iso\Ext^1(\sheaf{P}_{\pm a_i},\sheaf{P}_{\pm a_i})
\into \Ext^1(\sheaf{B}, \sheaf{B})
\end{align*}
induced by projection to and inclusion of the summand $\sheaf{P}_{\pm a_i}$ of $\sheaf{B}$.
Note that also the direct sum of the $g_i^\pm$'s is an isomorphism, since 
$\Ext^1(\sheaf{P}_x,\sheaf{P}_y) = H^1(X, \sheaf{P}_{y-x})$ vanishes unless $x=y$.

The obstruction map $\ob$ takes values in homotopy classes of morphisms
$M^\spot\to M^\spot[2]$ of complexes, modulo homotopy. Such a morphism
is given by a single homomorphism
from $\sheaf{A}$ to $\sheaf{C}$, which can be presented as an $(N-1)\times(N-1)$ matrix
with entries in $\Gamma(X,\OO_X(2\Theta))$. We now give such a matrix
representative for the homotopy class $\ob(\xi)$, for any element $\xi$
in each summand $H^1(X,\OO_X)$ of $E^{0,1}_2$.

\begin{lemma}\label{lemma:H1}
For every $i$, the boundary map of the long exact cohomology sequence
associated to the Koszul complex
\begin{equation*}
\begin{diagram}
0 &
\rTo & \OO_X &
\rTo^\koszA & \OO_X(\Theta_{a_i})\plus\OO_X(\Theta_{-a_i}) &
\rTo^\koszB & \sheaf{I}_{Y_i}(2\Theta) &
\rTo & 0
\end{diagram}
\end{equation*}
induces an isomorphism
\begin{equation*}
H^0(X, \sheaf{I}_{Y_i}(2\Theta))/\langle\vartheta^+_i\vartheta^-_i\rangle \iso
H^1(X, \OO_X)
\end{equation*}
where $\langle\vartheta^+_i\vartheta^-_i\rangle$ denotes the one dimensional vector
space spanned by the section $\vartheta^+_i\vartheta^-_i$.
\end{lemma}

\begin{proof}
Since $H^1(X,\OO_X(\Theta_{\pm a_i}))=0$, there is an induced right exact sequence
\begin{equation*}
\begin{diagram}
\begin{matrix}
H^0(X,\OO_X(\Theta_{a_i}))\\
\plus\\
H^0(X,\OO_X(\Theta_{-a_i}))
\end{matrix}
& \rTo &
H^0(X,\sheaf{I}_{Y_i}(2\Theta))
& \rTo &
H^1(X,\OO_X)
& \rTo & 0.
\end{diagram}
\end{equation*}
Each summand $H^0(X,\OO_X(\Theta_{\pm a_i}))$ is spanned by
$\vartheta^{\pm}_i$, which is sent
to $\mp\vartheta^+_i\vartheta^-_i$ in $H^0(X,\sheaf{I}_{Y_i}(2\Theta))$.
\end{proof}

\begin{proposition}\label{prop:ob}
Let $\xi\in H^1(X,\OO_X)$.
The obstruction map $\ob$ does the following on each summand in its domain:
\begin{enumerate}
\item Lift $\xi$ to sections $u$ and $v$
of $\sheaf{I}_{Y_i}(2\Theta)$ and $\sheaf{I}_{Y_N}(2\Theta)$, respectively, using
the lemma. Then
$\ob(f_{ij}(\xi))$ is represented by the $(N-1)\times(N-1)$ matrix having $j$'th column
(the transpose of)
\begin{equation*}
\begin{matrix}
(v & \cdots & v & u+v & v & \cdots & v)\\
    &        &    & \stackrel{\uparrow}{\scriptstyle\text{(entry $i$)}}
\end{matrix}
\end{equation*}
and zeros everywhere else.
\item Lift $\xi$ to a section $u$ of $\sheaf{I}_{Y_i}(2\Theta)$. If $i\ne N$, then
$\ob(g_i^\pm(\xi))$ is represented by the $(N-1)\times(N-1)$ matrix having $u$ at
entry $(i,i)$, and zeros everywhere else. The remaining case $\ob(g_N^\pm(\xi))$ is
represented by the $(N-1)\times(N-1)$ matrix having
all entries equal to $u$.
\item Lift $\xi$ to sections $u$ and $v$ of $\sheaf{I}_{Y_j}(2\Theta)$ and
$\sheaf{I}_{Y_N}(2\Theta)$, respectively. Then $\ob(h_{ij}(\xi))$ is represented
by the $(N-1)\times(N-1)$ matrix having $i$'th row
\begin{equation*}
\begin{matrix}
(v & \cdots & v & u+v & v & \cdots & v)\\
    &        &    & \stackrel{\uparrow}{\scriptstyle\text{(entry $j$)}}
\end{matrix}
\end{equation*}
and zeros everywhere else.
\end{enumerate}
\end{proposition}

\begin{proof}
Notation: In the commutative diagrams that follow, we will use dotted
arrows roughly to indicate maps that are not given to us, but need to be filled in
by some construction.

\textbf{Part 1:} We view $\xi$ as an extension in $\Ext^1(\OO_X(-\Theta),\OO_X(-\Theta))$.
Writing out the description of $\ob$ from Lemma \ref{lemma:obstruction}
in this situation, one arrives at the diagram
\begin{equation}\label{eq:extension1}
\begin{diagram}[width=3em]
0 & \rTo & \OO_X(-\Theta) & \rTo           & \sheaf{F}           & \rTo & \OO_X(-\Theta) & \rTo & 0\\
  &      & \dInto         & \rdTo^{\phi_i} & \dDashto_{\hat{\phi}_i} &      & \dDashto\\
  &      & \sheaf{A}              & \rTo_\phi      & \sheaf{B}                   & \rTo_{\psi} & \sheaf{C}
\end{diagram}
\end{equation}
constructed as follows: The top row is the extension $\xi$ and the bottom row is the monad.
The leftmost vertical map is inclusion of the $i$'th summand, so in terms of matrices,
$\phi_i$ is
the $i$'th column of $\phi$. We are required to extend $\phi_i$ to a map $\hat{\phi}_i$
making the left part of the diagram commute: one way of doing this is detailed below.
The induced vertical map on the right, precomposed with projection $\sheaf{A}\to \OO_X(-\Theta)$
on the $j$'th summand, is a representative for $\ob(f_{ij}(\xi))$.

The assumption that $u$ is a lifting of $\xi$, means that there is a commutative
diagram
\begin{equation*}
\begin{diagram}
0 & \rTo & \OO_X(-\Theta) &
\rTo^{\left(
\begin{smallmatrix}
\vartheta_i^+\\
\vartheta_i^-\\
\end{smallmatrix}
\right)}
& \sheaf{P}_{a_i}\plus \sheaf{P}_{-a_i} &
\rTo^{\left(
\begin{smallmatrix}
\vartheta_i^- & -\vartheta_i^+
\end{smallmatrix}
\right)}
& \sheaf{I}_{Y_i}(\Theta) & \rTo 0 \\
  &      & \uEq           &      & \uTo^{\tilde{u}}      &      & \uTo^u \\
0 & \rTo & \OO_X(-\Theta) & \rTo & \sheaf{F}             & \rTo & \OO_X(-\Theta) & \rTo 0
\end{diagram}
\end{equation*}
in which the rightmost square is a pullback. Similarly, the section $v$
fits in the pullback diagram:
\begin{equation*}
\begin{diagram}
0 & \rTo & \OO_X(-\Theta) &
\rTo^{\left(
\begin{smallmatrix}
\vartheta_N^+\\
\vartheta_N^-\\
\end{smallmatrix}
\right)}
& \sheaf{P}_{a_N}\plus \sheaf{P}_{-a_N} &
\rTo^{\left(
\begin{smallmatrix}
\vartheta_N^- & -\vartheta_N^+
\end{smallmatrix}
\right)}
& \sheaf{I}_{Y_N}(\Theta) & \rTo 0 \\
  &      & \uEq           &      & \uTo^{\tilde{v}}      &      & \uTo^v \\
0 & \rTo & \OO_X(-\Theta) & \rTo & \sheaf{F}             & \rTo & \OO_X(-\Theta) & \rTo 0
\end{diagram}
\end{equation*}
Now define $\hat{\phi}_i$ to be
\begin{equation*}
\begin{diagram}
(\tilde{u},\tilde{v})\colon & \sheaf{F} & \rTo & (\sheaf{P}_{a_i}\plus \sheaf{P}_{-a_i})\plus (\sheaf{P}_{a_N}\plus \sheaf{P}_{-a_N})
\end{diagram}
\end{equation*}
followed by the appropriate inclusion to $\sheaf{B}$. One verifies immediately
that $\hat{\phi}_i$ extends $\phi_i$ in \eqref{eq:extension1}, and that
the induced map in the rightmost part of that diagram is given by the
vector as claimed in part 1.

\textbf{Part 2:} We view $\xi$ as an extension in $\Ext^1(\sheaf{P}_{\pm a_i},\sheaf{P}_{\pm a_i})$.
In this situation, the description of $\ob$ from Lemma
\ref{lemma:obstruction} boils down to a diagram
\begin{equation*}
\begin{diagram}[width=3em]
0 & \rTo & \sheaf{P}_{\pm a_i} & \rTo & \sheaf{G} & \rTo & \sheaf{P}_{\pm a_i} & \rTo & 0\\
  &      & \dTo^{\psi_i^\pm} &  \ldDashto  &           &   \luDashto   & \uTo_{\phi_i^\pm}\\
  &      & \sheaf{C} & &  & & \sheaf{A}
\end{diagram}
\end{equation*}
constructed as follows: The top row is $\xi$ and the
two vertical maps $\phi_i^\pm$ and $\psi_i^\pm$ denote the $i$'th row of $\phi^\pm$
and the $i$'th column of $\psi^\pm$. The task is to lift $\phi_i^\pm$ to the
rightmost dotted arrow, and to extend $\psi_i^\pm$ to the leftmost dotted arrow.
The composition of the two dotted arrows is then a representative for
$\ob(g_i^\pm(\xi)))$.

First consider the problem of lifting $\mp\vartheta_i^\pm$ and extending $\vartheta_i^\mp$
to maps $s$ and $t$ as in the following diagram:
\begin{equation}\label{eq:st-diagram}
\begin{diagram}[width=3em]
0 & \rTo & \sheaf{P}_{\pm a_i} & \rTo & \sheaf{G} & \rTo & \sheaf{P}_{\pm a_i} & \rTo & 0\\
  &      & \dTo^{\vartheta_i^\mp} &  \ldDashto_t  &    &   \luDashto_s   & \uTo_{\mp\vartheta_i^\pm}\\
  &      & \OO_X(-\Theta) &  &  &  & \OO_X(\Theta)
\end{diagram}
\end{equation}
Suppose such a diagram is given. If $i<N$, then
$\mp(0,\dots,0,s,0,\dots,0)$ would lift $\phi_i^\pm$ and $\pm(0,\dots,0,t,0,\dots,0)$
would extend $\psi_i^\pm$. Their composition is the matrix having $t\circ s$ in
entry $(i,i)$, and zeros elsewhere. If $i=N$, then similarly $\mp(s,\dots,s)$
and $\pm(t,\dots,t)$ would be the required lift and extension. Their composition
is the matrix having all entries equal to $t\circ s$. Thus part 2 of the
proposition will be established once we have constructed such maps $s$ and $t$
having composition $t\circ s = u$.

Now use that $\xi$ is the pullback of the Koszul complex for $Y_i$ along $u$.
This enables us to construct the commutative diagram
\begin{equation*}
\newcommand{\kzA}{\left(
\begin{smallmatrix}
\vartheta_i^+\\
\vartheta_i^-
\end{smallmatrix}
\right)}
\newcommand{\kzB}{\left(
\begin{smallmatrix}
\vartheta_i^- & -\vartheta_i^+
\end{smallmatrix}
\right)}
\begin{diagram}[tight,width=2.5em,height=1.75em]
0 & \rTo & \OO_X(-\Theta) && \rTo^\kzA && \sheaf{P}_{a_i}\plus \sheaf{P}_{-a_i} && \rTo^\kzB && \sheaf{I}_{Y_i}(\Theta) & \rTo & 0 \\
  &      & \dEq           &  \rdTo_{\vartheta_i^\mp} & & \ldTo && \luTo & & \ruTo_{\mp\vartheta_i^\pm}     & \uTo_u\\
  &      &                &  & \sheaf{P}_{\mp a_i}  &&\uTo&&      \sheaf{P}_{\mp a_i}\\
  &      &                && \dEq     &&            &  &\uTo^u &&\\
0 & \rTo & \OO_X(-\Theta) & \rLine & \VonH & \rTo & \sheaf{G}' & \rLine & \VonH & \rTo &  \OO_X(-\Theta) & \rTo & 0 \\
  &      &                &  \rdTo_{\vartheta_i^\mp} & & \ldDashto && \luDashto &  & \ruTo_{\mp\vartheta_i^\pm}      \\
  &      &                &  & \sheaf{P}_{\mp a_i}  &&&&      \sheaf{P}_{\mp a_i}(-2\Theta)
\end{diagram}
\end{equation*}
as follows: The top row is the Koszul complex. In the top part, the
unlabelled diagonal arrows are the canonical inclusion of and projection to the
summand $\sheaf{P}_{\mp a_i}$. In particular their composition is the identity
map. Pull back along $u$ to get the short exact sequence in the lower
part of the diagram. Thus this sequence coincides with $\xi$, twisted
by $\sheaf{P}_{\mp a_i}(-\Theta)$. There are now uniquely determined dotted
arrows making the diagram commute, and their composition is $u$. Twisting
back by $\sheaf{P}_{\pm a_i}(\Theta)$, the lower part of the
diagram is thus the required diagram \eqref{eq:st-diagram}. This ends the proof of part 2.

\textbf{Part 3} is essentially dual to part 1, and is left out.
\end{proof}

By Lemma \ref{lemma:theta-intersect}, we have
\begin{equation}\label{eq:sum}
H^0(\sheaf{I}_{Y_i}(2\Theta))\plus H^0(\sheaf{I}_{Y_j}(2\Theta)) = H^0(\OO_X(2\Theta))
\end{equation}
for all $i\ne j$ (the Lemma gives an inclusion of the left hand side into the
right hand side, and by Riemann-Roch, the two sides have the same dimension). This decomposition of sections of $\OO_X(2\Theta)$, together
with the explicit description of $\ob$ in the proposition, enables us to conclude:

\begin{corollary}
The obstruction map $\ob$ is surjective.
\end{corollary}

\begin{proof}
We show that any $(N-1)\times(N-1)$
matrix of sections of $\OO_X(2\Theta)$ represents an element in
the image of $\ob$. Let $i$ and $j$ be arbitrary indices between $1$ and $N-1$.

Step 1: Let $v$ be a section of $\sheaf{I}_{Y_N}(2\Theta)$.
Let $\xi\in H^1(X,\OO_X)$ be the image of $v$ under the boundary map in Lemma
\ref{lemma:H1}, and lift
$\xi$ to another section $u$ of $\sheaf{I}_{Y_i}(2\Theta)$. By parts 1 and 2 of
the Proposition, 
$\ob(f_{ii}(\xi)-g_i^\pm(\xi))$ is represented by the matrix having zeros except
for in column $i$, where all elements equal $v$. Similarly
$\ob(h_{ii}(\xi)-g_i^\pm(\xi))$ is represented by a matrix
having all entries of row $i$ equal to $v$, and zeros elsewhere.

Step 2: Let $u$ be a section of $\sheaf{I}_{Y_i}(2\Theta)$. By part 1 of the Proposition and the
previous step, we can find a matrix representing an element in the image of
$\ob$, with $u$ as entry $(i,j)$ and
zeros elsewhere. Similarly, for any section $u$ of $\sheaf{I}_{Y_j}(2\Theta)$,
combining part 3 of the proposition with the previous step, we obtain a matrix having $u$ as entry
$(i,j)$ and zeros elsewhere.

Step 3: If $i\ne j$, the previous step and \eqref{eq:sum}
enables us to construct a matrix with arbitrary
entries outside the diagonal. Combining this with step 1, we can construct a matrix
having any given section of $\sheaf{I}_{Y_N}(2\Theta)$ at entry $(i,i)$,
and zeros elsewhere.
By step 2 we can also construct a matrix having any given section of $\sheaf{I}_{Y_i}(2\Theta)$
at $(i,i)$, and zeros elsewhere. By \eqref{eq:sum} with $j=N$,
this enables us to hit arbitrary elements along the diagonal, too.
\end{proof}

\begin{corollary}
The spectral sequence \eqref{eq:spectral} degenerates at $E_3$.
\end{corollary}

\begin{proof}
The previous corollary implies $E^{2,0}_3=0$. By duality also
$E^{-2,3}_3=0$. It follows from the shape of the first sheet,
Figure \ref{fig:spectral}, that all differentials vanish at the
$E_3$-level and beyond.
\end{proof}

\subsection{First order deformations}

From the calculations in the previous section, we can understand infinitesimal
deformations of the vector bundle $\sheaf{E}$ in terms of its monad. Let
$k[\epsilon]$ be the ring of dual numbers. By a first order deformation of
$M^\spot$, we mean a monad over $X\tensor_k k[\epsilon]$, with $M^\spot$ as
fibre over $\epsilon=0$, modulo isomorphism.

\begin{theorem}\label{thm:qualitative}
Let $M^\spot$ be a decomposable monad with cohomology $\sheaf{E}$.
The vector spaces of first order infinitesimal deformations
of $M^\spot$ and of $\sheaf{E}$ are isomorphic via the
natural map, sending a first order deformation of $M^\spot$ to its
cohomology.
\end{theorem}

\begin{proof}
Since the hyperext spectral sequence associated to
the monad degenerates at $E_3$,
and the only $E^{pq}_3$ terms
with $p+q=1$ are $E^{0,1}_3$ and $E^{1,0}_3$, there is a short
exact sequence
\begin{equation}\label{eq:ses}
\begin{diagram}
0 & \rTo & E^{1,0}_3 & \rTo & \Ext^1(\sheaf{E},\sheaf{E}) & \rTo & E^{0,1}_3 & \rTo & 0.
\end{diagram}
\end{equation}
Let $D(M^\spot)$ be the vector space of first order
deformations of $M^\spot$. Thus the claim is
that the natural map $D(M^\spot)\to \Ext^1(\sheaf{E},\sheaf{E})$
is an isomorphism. It suffices to show that $D(M^\spot)\to E^{0,1}_3$
is surjective, and that its kernel maps isomorphically to $E^{1,0}_3$.

Now $E^{0,1}_3$ is the kernel of the obstruction map $\ob=d^{0,1}_2$.
By Lemma \ref{lemma:obstruction}, this is the space of those first order
deformations of the objects in $M^\spot$, that allow the differential
$d_M$ to extend (non uniquely) to the deformed objects. Via this
identification, $D(M^\spot) \to E^{0,1}_3$ is the natural forgetful
map, so it is surjective. Moreover, its kernel is the space of first
order deformations of the differential in $M^\spot$, keeping the
objects fixed. It remains to see that this space gets identified
with $E^{1,0}_3$.

By the shape of the spectral sequence (Figure \ref{fig:spectral}) we have $E^{1,0}_3 = E^{1,0}_2$,
and, by Lemma \ref{lem:htpy}, this is
\begin{equation*}
E^{1,0}_2 = \Hom_{K(X)}(M^\spot,M^\spot[1]).
\end{equation*}
The inclusion of $E^{1,0}_2$ into $\Ext^1(\sheaf{E},\sheaf{E})$
is the 
edge map discussed in Section \ref{sec:edge}, i.e.~the canonical
map
\begin{equation}\label{eq:edge}
\Hom_{K(X)}(M^\spot, M^\spot[1]) \to \Hom_{D(X)}(M^\spot, M^\spot[1]).
\end{equation}
This can be factored as follows: a morphism of complexes in an arbitrary
abelian category $f\colon X^\spot \to Y^\spot[1]$
gives rise to a short exact sequence of complexes
\begin{equation}\label{eq:cone-ses}
\begin{diagram}
0 & \rTo & Y^\spot & \rTo^{\beta} & Z^\spot & \rTo^{\alpha} & X^\spot & \rTo & 0
\end{diagram}
\end{equation}
where $Z^\spot = C(f[-1])$ is the mapping cone of $f[-1]$, which has objects $X^i\plus
Y^i$ in degree $i$ and differential $(x,y) \mapsto (dx, f(x)+dy)$. 
The maps $\alpha$ and $\beta$ are the canonical ones.
Moreover, the usual Yoneda construction of elements in $\Ext^1$
from short exact sequences (of objects) can be extended to complexes,
by associating to any short sequence \eqref{eq:cone-ses} the roof
\begin{equation*}
\begin{diagram}
& & C(\beta) \\
& \ldTo^{\text{qism}} & & \rdTo \\
X^\spot & & & & Y^\spot[1]
\end{diagram}
\end{equation*}
where $C(\beta)$ is the mapping cone, with objects $Y^{i+1}\plus Z^i$ in degree $i$
and differential $(y,z)\mapsto (-dy, \beta(y)+dz)$. The leftmost map is given
by projection, and is a quasi-isomorphism, whereas the rightmost map is projection
followed by $\alpha$. This roof defines a morphism $X^\spot\to Y^\spot[1]$ in the
derived category. Moreover, the diagram obtained from the roof by adding the
negative of the map $f\colon X^\spot\to Y^\spot[1]$ we started with, is commutative
up to homotopy, so the roof and $-f$ defines the same map in the derived category.

Thus we have factored the edge map \eqref{eq:edge} via short exact
sequences, by sending $f\colon M^\spot\to M^\spot[1]$ to
the short exact sequence
%
%
\begin{equation*}
\begin{diagram}
0 & \rTo & M^\spot & \rTo & C(f[-1]) & \rTo & M^\spot & \rTo & 0.
\end{diagram}
\end{equation*}
The associated element in $\Hom_{D(X)}(M^\spot, M^\spot[1])$ corresponds, up to sign,
to the Yoneda class in $\Ext^1(\sheaf{E},\sheaf{E})$ obtained by taking
the $H^0$ cohomology of each complex in this short exact sequence. To phrase
this in terms of first order deformations, we 
rewrite the cone
$C(f[-1])$ as the complex $M^\spot\tensor_k k[\epsilon]$ equipped with the
differential $d_M\tensor 1 + f\tensor \epsilon$. The corresponding deformation
of $\sheaf{E}$ is the $H^0$ cohomology of this complex. But this is the
required result, since any
differential on $M^\spot\tensor_k k[\epsilon]$ that specializes to $d_M$
for $\epsilon = 0$ has the form $d_M\tensor 1 + f\tensor \epsilon$, for some $f$ satisfying
\begin{equation*}
(d_M\tensor 1 + f\tensor\epsilon)^2 = 0.
\end{equation*}
Since $d_M^2=0$ and $\epsilon^2=0$, this says that $f\circ d_M + d_M\circ f = 0$,
i.e.~$f$ defines a morphism $M^\spot\to M^\spot[1]$. This gives the required
identification between $E^{1,0}_2$ and deformations of the differential.
\end{proof}

Next, we give the dimension formula for $\Ext^1(\sheaf{E},\sheaf{E})$,
which we phrase in a twist invariant way.

\begin{theorem}\label{thm:quantitative}
Let $\sheaf{E}$ be a rank $2$ vector bundle obtained as
the cohomology of a decomposable monad, or the twist of
such a bundle by a line bundle. Then
\begin{equation*}
\dim \Ext^1(\sheaf{E},\sheaf{E}) =
\tfrac{1}{3}\Delta(\sheaf{E})\cdot\Theta + 5
\end{equation*}
where $\Delta$ denotes the discriminant $4c_2-c_1^2$.
\end{theorem}

\begin{proof}
Both sides of the equation are invariant under twist, so it suffices to
verify the formula when $\sheaf{E}$ is the cohomology of a decomposable
monad. Consider again the short exact sequence \eqref{eq:ses}.

The space $E^{0,1}_3$ is the kernel of the map $\ob = d^{0,1}_2$ studied
in Section \ref{sec:E2diff}. Its domain \eqref{eq:E012} has dimension
\begin{equation*}
(2(N-1)^2 + 2N) \dim H^1(X, \OO_X) = 6(N-1)^2 + 6N
\end{equation*}
and its codomain $E^{2,0}_2$ has dimension $6(N-1)^2 - N + 2$, by Lemma \ref{lem:htpy}.
Since $\ob$ is surjective, the dimension of its kernel $E^{0,1}_3$ is thus $7N-2$.
Moreover, the dimension of $E^{1,0}_3=E^{1,0}_2$ is $N-1$ by the same Lemma, so
$\Ext^1(\sheaf{E},\sheaf{E})$ has dimension $8N-3$.

On the other hand,
we know from the Serre construction that $\sheaf{E}(\Theta)$ has Chern
classes $c_1=2\Theta$ and $c_2=N\Theta^2$, and thus discriminant
$(4N-4)\Theta^2$. The formula follows.
\end{proof}

\begin{remark}
The space of first order deformations obtained by varying the
isomorphism $\omega_Y\iso \OO_Y(2\Theta)$, coincides with the space
of first order deformations of the differential in $M^\spot$.
In fact, it is trivial that the former is contained in the latter,
and these spaces have the same dimension $N-1$, using Proposition
\ref{prop:serrefamily}.
\end{remark}

\begin{remark}
The short exact sequence \eqref{eq:ses}, and its interpretation given
in the proof of Theorem \ref{thm:qualitative},
is not intrinsic to $\sheaf{E}$, but results from our choice of representing
$\sheaf{E}$ by a decomposable monad.
However, deformation of the differential in
$M^\spot$, or
equivalently, variation of the isomorphism $\OO_Y(2\Theta)\iso \omega_Y$
in the Serre construction, defines a rational $(N-1)$-dimensional
subvariety through $\sheaf{E}$ in its moduli space, whose tangent space
is $E^{0,1}_3$. It seems plausible that this $(N-1)$-dimensional
rational variety can be intrinsically characterized as
the unique (maximal) rational variety through $\sheaf{E}$.
\end{remark}

\section{Birational description of $M(0,\Theta^2)$}

As before, let $(X,\Theta)$ be a principally polarized abelian threefold
with Picard number $1$.
We write $M(c_1,c_2)$ for the coarse moduli space of
stable rank $2$ vector bundles on $X$ with the indicated Chern classes.

The main point in the preceding section is that all first order
infinitesimal deformations of the vector bundles constructed in Section
\ref{sec:duality}, in the case of even
$c_1$, can be realized as first order infinitesimal deformations of a
monad. In this section we show that in the first nontrivial example,
corresponding to $N=2$, this statement holds not only infinitesimally,
but Zariski locally: by deforming the monad, we realize a Zariski open
neighbourhood of the vector bundle in its moduli space. In terms of
the Serre construction, this is the case corresponding to curves $Y_1\union Y_2$
with two components $Y_i = \Theta_{a_i}\intersect\Theta_{-a_i}$.

\begin{theorem}
Let $\sheaf{E}$ be the rank $2$ cohomology vector bundle 
of a decomposable monad, as in Proposition \ref{prop:decomposable} for $N=2$.
Then, Zariski locally around $\sheaf{E}$, the moduli space
$M(0,\Theta^2)$ is a uniruled, nonsingular variety of dimension $13$.

More precisely, there is a Zariski open neighbourhood around $\sheaf{E}$
which is isomorphic to a nonsingular Zariski open subset of a
$\PP^1$-bundle over a finite quotient of $X^2\times_X X^2\times_X X^2$,
where $X^2$ is considered as a scheme over $X$ via the group law.
\end{theorem}

\begin{proof}
We write down a parameter space for the family of monads
\begin{equation}\label{eq:N2monad}
\begin{diagram}
0 & \rTo & \sheaf{P}_{b'}(-\Theta) & \rTo^\phi & \Plus_{i=1}^2 (\sheaf{P}_{a_i}\plus \sheaf{P}_{a_i'}) & \rTo^\psi & \sheaf{P}_{b}(\Theta) & \rTo & 0
\end{diagram}
\end{equation}
where $a_i, a_i', b, b'$ are sufficiently general points in $X$ satisfying
\begin{equation}\label{eq:N2equations}
a_1+a_1' = a_2+a_2' = b+b',
\end{equation}
and
\begin{equation*}
\phi = (\vartheta_1, \vartheta_1', \vartheta_2, \vartheta_2'),\quad
\psi = (\vartheta_1', -\vartheta_1, \vartheta_2', -\vartheta_2)
\end{equation*}
and where the $\vartheta$'s are required to be nonzero, but otherwise
arbitrary.

Viewing $X^2$ as a variety over $X$ via the group law, the fibred product
$X^2\times_X X^2\times_X X^2$ is the subvariety of $X^6$ defined by
\eqref{eq:N2equations}. Let
\begin{equation*}
T\subset X^2\times_X X^2\times_X X^2
\end{equation*}
be the open subset consisting of sixtuples $(a_1,a_1';a_2,a_2';b,b')$ where
the leading four entries are all distinct. Later we may have to shrink $T$ further.
With the help of the Poincar\'e line bundle on $X\times X$ it is clear that,
on $T\times X$, there exist vector bundles $\sheaf{A}$, $\sheaf{B}$, $\sheaf{C}$
whose fibres over a sixtuple in $T$ are the three objects in \eqref{eq:N2monad}.
The sixtuples $(a_1,a_1';a_2,a_2;b,b')$ in $T$ corresponding to the same three objects
constitute an orbit for the action of
\begin{equation}\label{eq:Ggroup}
G = (\ZZ/(2)\plus\ZZ/(2)) \rtimes \ZZ/(2)
\end{equation}
on $T$, where the action of the first semidirect factor is given by the transpositions
$a_1\leftrightarrow a_1'$ and $a_2\leftrightarrow a_2'$, and the last factor
acts by $(a_1,a_1')\leftrightarrow (a_2,a_2')$. Thus $T/G$ is a parameter space for
the objects in \eqref{eq:N2monad}.

Next we parametrize the maps $\phi$ and $\psi$, which are given by four nonzero sections
\begin{equation}\label{eq:thetas}
\begin{aligned}
\vartheta_1 &\in \Gamma(X, \sheaf{P}_{a_1-b'}(\Theta)) &\vartheta'_1 &\in \Gamma(X, \sheaf{P}_{a'_1-b'}(\Theta))\\
\vartheta_2 &\in \Gamma(X, \sheaf{P}_{a_2-b'}(\Theta)) &\vartheta'_2 &\in \Gamma(X, \sheaf{P}_{a'_2-b'}(\Theta))
\end{aligned}
\end{equation}
There exist line bundles $L_1,L_1',L_2,L_2'$ on $T$, whose fibres over a sixtuple $(a_1,a_1';a_2,a_2';b,b')$
are these (one dimensional) spaces of global sections. Thus, writing
\begin{equation*}
F = \Plus_{i=1}^2 (L_i\plus L_i') \xrightarrow{p} T,
\end{equation*}
a point of $F$, whose four entries are all nonzero, corresponds to a monad \eqref{eq:N2monad}.
More precisely, writing $p_X$ for the product
$p\times\id_X\colon F\times X\to T\times X$, there exists a monad
\begin{equation*}
\begin{diagram}
p_X^*\sheaf{A} & \rTo^\Phi & p_X^*\sheaf{B} & \rTo^\Psi & p_X^*\sheaf{C}
\end{diagram}
\end{equation*}
on $F\times X$, whose restriction to the point in $F$ given by \eqref{eq:thetas} is
\eqref{eq:N2monad}. Let $F'\subset F$ be the open subset consisting of quadruples with only
nonzero entries. The cohomology of the ``universal'' monad above is a family of vector bundles over $F'$,
giving rise to a morphism of schemes
\begin{equation}\label{eq:maptoM}
\phi\colon F' \to M(0,\Theta^2).
\end{equation}
To make this morphism an embedding, we will divide by the group $G$ to get rid of the ambiguity
in the parametrization of the objects by $T$, and further divide by another group $\Gamma$ to take
care of distinct maps $\phi$, $\psi$ which give isomorphic monads.

For a fixed base point in $T$, and hence fixed objects in \eqref{eq:N2monad},
the tuples $(\vartheta_1,\vartheta_1',\vartheta_2,\vartheta_2')$ which define isomorphic monads
constitute orbits under the
following group action on $F$: view $\Gm^2$ as a variety over $\Gm$ via the multiplication map, and
let
\begin{equation*}
\Gamma=\Gm^2\times_{\Gm} \Gm^2.
\end{equation*}
Its closed points are tuples
$(\lambda_1,\lambda_1';\lambda_2,\lambda_2')$ satisfying
$\lambda_1\lambda_1'=\lambda_2\lambda_2'$. The action on the fibres of $F$ is
given on closed points by
\begin{equation*}
(\vartheta_1,\vartheta_1',\vartheta_2,\vartheta_2')
\mapsto (\lambda_1\vartheta_1,\lambda_1'\vartheta_1',\lambda_2\vartheta_2,\lambda_2'\vartheta_2').
\end{equation*}
There is a short exact sequence of group varieties
\begin{equation*}
\begin{diagram}
1 & \rTo & \Gm^2 & \rTo & \Gamma & \rTo & \Gm & \rTo 1
\end{diagram}
\end{equation*}
where the inclusion sends $(\lambda_1,\lambda_1')$ to $(\lambda_1,\lambda_1^{-1},\lambda_2,\lambda_2^{-1})$
and the projection sends $(\lambda_1,\lambda_1',\lambda_2,\lambda_2')$ to $\lambda_i\lambda_i'$.
Correspondingly, we determine $F/\Gamma$ in two steps. Firstly, the categorical quotient by
the $\Gm^2$-action is
\begin{equation*}
F/\Gm^2 \iso \Plus_{i=1}^2 (L_i\tensor L_i')
\end{equation*}
and the quotient map $F\to F/\Gm^2$ corresponds to multiplication in the fibres (locally on $T$,
this is the product of two copies of the quotient $\Spec R[x,y] \to \Spec R[xy]$
for the $\Gm$-action $(x,y) \to (\lambda x,\lambda^{-1} x)$ on $\AA^2_R$ over an arbitrary ring $R$).
The induced action of $\Gamma/\Gm^2\iso\Gm$ on the rank two vector bundle $F/\Gm^2$ is
multiplication in the fibres, so the quotient $P=\{F\setminus 0\}/\Gamma$ is
\begin{equation*}
P=\PP(\dual{\Plus_{i=1}^2 (L_i\tensor L_i')}),
\end{equation*}
which is a $\PP^1$-bundle over $T$. The image $P'\subset P$
of $F'\subset F$ is an open subset; in fact it is the complement of the two natural sections
corresponding to the subbundles $L_i\tensor L_i'$ of $F/\Gm^2$. The restricted quotient map
\begin{equation*}
F' \to F'/\Gamma = P'
\end{equation*}
is a geometric quotient; in particular its fibres are orbits in $F'$. It is clear that the morphism
\eqref{eq:maptoM} is invariant with respect to the $\Gamma$-action on $F'$, so there is an induced
morphism
\begin{equation*}
\overline{\phi}\colon P'\to M(0,\Theta^2).
\end{equation*}
Moreover, the (free) action \eqref{eq:Ggroup} of $G$ on $T$ has a canonical lift to $P$, and $P'$ is
$G$-invariant. Again $\overline{\phi}$ is invariant under this action, so we obtain the
$\PP^1$-bundle $P/G$ over $T/G$, together with an open subset $P'/G$ and an induced morphism
\begin{equation*}
\overline{\overline{\phi}}\colon P'/G \to M(0,\Theta^2).
\end{equation*}
By construction, the domain $P'/G$ parametrizes isomorphism classes of monads of the form
\eqref{eq:N2monad}. Given two such monads $M_1^\spot$ and $M_2^\spot$, with cohomology
$\sheaf{E}_1$ and $\sheaf{E}_2$, the first hyperext
spectral sequence gives an isomorphism $\Hom(M_1^\spot,M_2^\spot) \isoto \Hom(\sheaf{E}_1,\sheaf{E}_2)$.
Here, the domain is the $E^{0,0}_2$-term in the spectral sequence, which is
the group of morphisms of complexes (there are no homotopies, since $E^{-1,0}_1$ vanishes).
It follows
that $M^\spot_1$ and $M^\spot_2$ are isomorphic as complexes if and only if
$\sheaf{E}_1$ and $\sheaf{E}_2$ are isomorphic vector bundles.
In other words $\overline{\overline{\phi}}$ is
injective on closed points. Shrinking $T$ if necessary, we may apply Theorem \ref{thm:qualitative} to see
that $\overline{\overline{\phi}}$ is \'etale at points where $b=b'=0$. By shrinking its domain
further if necessary, we can assume that it is \'etale everywhere.
An \'etale and injective morphism is an open embedding, so we are done.
\end{proof}

\bibliographystyle{amsplain}
\bibliography{all}

\end{document}